\documentclass[referee,sn-basic]{sn-jnl}

\jyear{2022}%

\theoremstyle{thmstyleone}%
\newtheorem{theorem}{Theorem}
\newtheorem{lemma}{Lemma}%

\theoremstyle{thmstyletwo}%

\theoremstyle{thmstylethree}%

\usepackage{bm}
\def\blue#1{{\color{black}#1}}

\newcommand{\Var}{\mathop{\mathrm{Var}}}

\newcommand{\m}{\mathop{\bm{m}}}

\newcommand{\Gam}{\mathop{\bm{\Gamma}}}
\newcommand{\alp}{\mathop{\bm{\alpha}_0}}
\newcommand{\alpp}{\mathop{\bm{\alpha}_1}}

\newcommand{\ep}{\mathop{\bm{\epsilon}}}

\newcommand{\hGam}{\mathop{\hat{\bm{\Gamma}}}}
\newcommand{\halp}{\mathop{\hat{\bm{\alpha}}_0}}

\newcommand{\hbe}{\mathop{\hat{\bm{\beta}}}}

\newcommand{\0}{\mathop{\bf{0}}}

\newcommand{\Sigep}{\mathop{\bm{\Sigma}_{\epsilon}}}

\newcommand{\overD}{\mathop{\overset{D}{\rightarrow}}}

\begin{document}

\title[Generalized Varying Coefficient Mediation Models]{Generalized Varying Coefficient Mediation Models}


\author*[1]{\fnm{Jingyuan} \sur{Liu}}\email{jingyuan@xmu.edu.cn}

\author[2]{\fnm{Yujie} \sur{Liao}}\email{yxl629@psu.edu}

\author[2]{\fnm{Runze} \sur{Li}}\email{rzli@psu.edu}

\affil*[1]{\orgdiv{MOE Key Laboratory of Econometrics, Department of Statistics and Data Science,
School of Economics, Wang Yanan Institute for Studies in Economics, Fujian Key Lab of Statistics},
\orgname{Xiamen University},  \orgaddress{\city{Xiamen},  \state{Fujian 361005},  \country{China}}}


\affil[2]{\orgdiv{Department of Statistics}, \orgname{The Pennsylvania State University},
\orgaddress{\city{University Park}, \state{PA 16802}, \country{USA}}}


\abstract{Motivated by an analysis of causal mechanism from economic stress to entrepreneurial withdrawals through depressed affect, we develop a two-layer generalized varying coefficient mediation model. This model captures the bridging effects of mediators that may vary with another variable, by treating them as smooth functions of this variable. It also allows various response types by introducing the generalized varying coefficient model in the first layer. The varying direct and indirect effects are estimated through spline expansion. The theoretical properties of the estimated direct and indirect coefficient functions, including estimation biases, asymptotic distributions, and so forth, are explored. Simulation studies validate the finite-sample performance of the proposed estimation method. A real data analysis based on the proposed model discovers some interesting behavioral economic phenomenon, that self-efficacy influences the deleterious impact of economic stress, both directly and indirectly through depressed affect, on business owners' withdrawl intentions.}



\maketitle

\section{Introduction}\label{sec1}

Mediation analysis, which could be firstly traced back to decades ago \citep{wright1918, wright1934}, has been well developed and widely used for understanding the relationship between exposure variables and outcomes through some other variables called mediators. It assumes that mediators serve as a bridge along the causal chain from exposures to outcomes. Both the direct effects of the exposures on the outcome and the indirect effects of them through mediators can be assessed.
Statistical estimation and inference literature has made tremendous progress toward mediation analysis, ranging from fitting a succession of regression models \citep{baron1986,mackinnon2002comparison} to structural equations model (SEM) based methods \citep{Ditlevsen2005}. Mediation models with multiple or even high-dimensional mediators are also analysed \citep{VanderWeele2014, Huang2014, Huang2016hypothesis, Sohn2019, Zhou2020}. In addition, various outcome types have been systematically studied. For instance, \cite{VanderWeele2010odds} modified odds ratios in mediation models for a dichotomous outcome. \cite{Valeri2014} extended it to mediators with measurement errors. \cite{Wang2011estimating, Lange2011direct, Luo2020high}, among others, investigated mediation analysis in survival models.

Meanwhile, originating from psychology \citep{baron1986, mackinnon2008introduction}, mediation analysis has broad application prospects, including economics, behavioral science, epidemiology, biomedical research, and so forth. Inspired by an empirical study of behavioral economics \citep{pollack2012moderating}, we are motivated in this paper to explore the relationship between economic stress, entrepreneurs' depressed affect and their intentions to withdraw from business. According to \cite{baron2008role} and \cite{pollack2012moderating}, understanding the feelings and emotions of entrepreneurs in response to economic stress is salient in understanding their cognitive and behavioral processes such as opportunity recognition and resource acquisition. Keeping a business afloat under economic stress is challenging, especially in difficult economic times. During this process, consistent with the learned helplessness theory, entrepreneurs are likely to develop feelings such as depression, helplessness and hopelessness due to job-related insecurity \citep{seligman1972learned}. Further, such emotions may place substantial burden on individuals and compel entrepreneurs to change their business strategies and revenue generation mechanism \citep{egan2009coping, latham2009contrasting}. Therefore, we are motivated to adopt mediation models to investigate the relationship between the economic stress and entrepreneurs' withdrawal intentions mediated through depressed affect.

\begin{figure}[h]
\centering
\includegraphics[width=11cm,height=7cm]{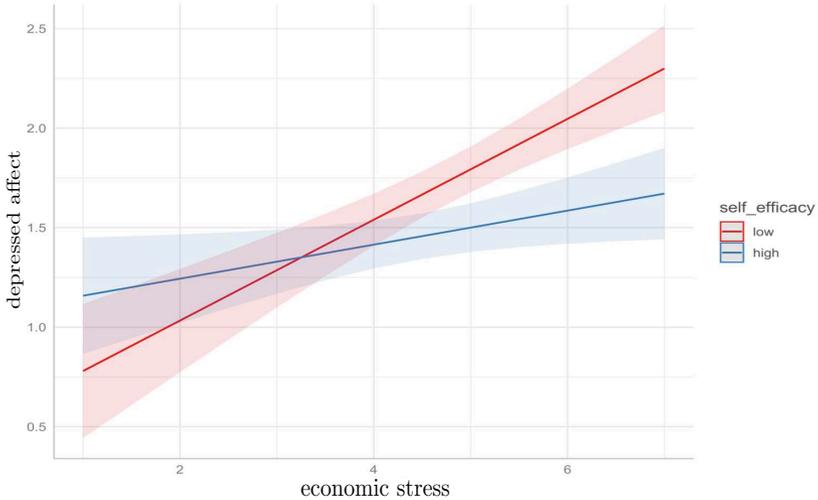}
\caption[Relationship between Economic Stress and Depressed Affect for Different Self-efficacy Levels]{Relationship between Economic Stress and Depressed Affect for Different Self-efficacy Levels}
\label{EDA_interaction}
\end{figure}

However, Figure \ref{EDA_interaction} shows that self-efficacy plays a crucial role in buffering the effect of economic stress on depressed affect. 
When stress becomes more severe, people with low self-efficacy are more likely to develop an increased depression. According to \cite{bandura1997self}, self-efficacy is a person's belief in their ability to achieve their goals successfully as a result of their own actions. It can influence people in several ways, such as thought patterns and how people deal with life stressors. Strong self-efficacy can provide critical confidence and alleviate feelings of depressed affect. People with a strong sense of self-efficacy are likely to master tasks and overcome abuse, while people without self-efficacy has little incentive to face difficulties  \citep{sequeira2007influence}. Therefore, the effect of economic stress on entrepreneurs' withdrawal intentions, both direct and indirect, might vary with self-efficacy of entrepreneurs. Some existing works, but not many, are applicable for the varying pattern. For instance, \cite{VanderWeele2017} tackled the time-varying exposures and mediators by a so-called ``mediational g-formula" using two marginal structural models. \cite{Liao2021varying} proposed using nonparametric procedures for estimating the varying direct and indirect effects when the outcome is continuous. However, neither of them is applicable for the current problem where the outcome withdrawal intention is dichotomous.

This motivates us to propose a two-layer generalized varying coefficient mediation model (GVCMM), in which the effect of exposures are represented as smooth functions of another effect modifier, while various types of outcomes can be accommodated. 
Not much work has been done in literature, to our best knowledge, to estimate the varying effects in generalized regression models. One main reference is  \cite{Verhasselt2014}, where they dealt with nonparametric smoothing and variable selection in single-layer generalized varying coefficient models. However, it is not applicable for mediation analysis. One main but not only obstacle is that in the two-layer mediation model, the mediators are random and hence typically unbounded, thus the existing theory in \cite{Verhasselt2014} with bounded covariates are not in point. For fitting GVCMM, polynomial spline methods \citep{huang2004polynomial, huang2002varying} are adopted to estimate the varying direct and indirect effects. The asymptotic normalities of estimated effects are also established. This builds up foundation for subsequential statistical inference of the varying direct and indirect effects.

The rest of the article is organized as follows. In Section 2, we propose a new generalized varying coefficient mediation model. We also develop an estimation procedure for indirect and direct effects, and establish asymptotic properties of the proposed estimators. In Section 3, simulation studies are conducted to assess the finite performance of the proposed procedure. In Section 4, the indirect effect of economic stress on entrepreneurial withdrawals via depressed affect, as well as its direct effect, are examined, and the buffering effect of self-efficacy is explored in details. Conclusion and discussion are given in Section 5. Technical details can be found in the appendix.

\section{Modeling Procedures}\label{sec2}
Let $y$ be the response, $\m$ be a vector of $p$ mediators, $\bm{x}$ be a vector of $q$ exposures, and $u$ be a univariate covariate. The generalized varying coefficient mediation model (GVCMM) assumes that
\begin{align}
\label{GVCMMmodel_y}g\{\mathrm{E}(y\mid u,\bm{m}, \bm{x})\}&=\bm{\alpha}_0(u)^T \bm{m} +\bm{ \alpha}_1(u)^T \bm{x}\\
\label{GVCMMmodel_m}\bm{m}&=\bm{\Gamma}(u)^T \bm{x} + \bm{\epsilon},
\end{align}
where $\ep$ is an error term with mean $\bm{0}$ and covariance matrix $\Sigep$, which is independent of $\m$, $\bm{x}$ and $u$. $g(\cdot)$ is a known link function - this broadens the application of the proposed mediation models, including but not limited to that with categorical or ordered response $y$. $\bm{\alpha}_0(\cdot)$, $\bm{\alpha}_1(\cdot)$ and $\bm{\Gamma}(\cdot)$ are unspecified smooth regression coefficient functions. Specifically, the GVCMM considers that conditional on $\bm{m}$, $\bm{x}$ and $u$, the distribution of $y$ belongs to an exponential family with density
$$c(y,\phi)\exp\{\frac{y\theta-b(\theta)}{\phi}\}.$$
With the canonical link,
$$\theta(u,\bm{m},\bm{x})=\bm{\alpha}_0(u)^T \bm{m} +\bm{ \alpha}_1(u)^T \bm{x}.$$

By the constructed models, $\bm{x}$ affects $y$ both directly through model  (\ref{GVCMMmodel_y}) and indirectly via the mediator vector $\m$ in model  (\ref{GVCMMmodel_m}). $\bm{\alpha}_1(u)$ is thus referred to as the direct effect of $\bm{x}$ on $y$, which may vary with $u$. To define the indirect effect, (\ref{GVCMMmodel_m}) indicates that
\begin{align}
\label{GVCMMmodel_indirect}\bm{\alpha}_0(u)^T\bm{m}&=\{\bm{\Gamma}(u)\bm{\alpha}_0(u)\}^T \bm{x} +\bm{ \alpha}_0(u)^T \bm{\epsilon}.
\end{align}
Hence we call $\bm{\beta}(u)=\bm{\Gamma}(u)\bm{\alpha}_0(u)$ the indirect effect of $\bm{x}$ on $y$.

\subsection{Estimation Procedures}

We in this section discuss the estimation procedures of the direct effect $\bm{\alpha}_1(u)$ and the indirect effect $\bm{\beta}(u)$. Notice that $\bm{\beta}(u)=\bm{\Gamma}(u)\bm{\alpha}_0(u)$, thus the estimation of $\bm{\beta}(u)$ essentially rely on that of $\bm{\Gamma}(u)$ and $\bm{\alpha}_0(u)$.

Suppose that $\{(u_i,\bm{m}_i,\bm{x}_i,y_i)\}_{i=1}^n$ is a random sample. We then rewrite model (\ref{GVCMMmodel_y}) and (\ref{GVCMMmodel_m}) in an elementwise manner. That is,
\begin{align*}
    g\{\mathrm{E}(y_i\mid u_i,\bm{m}_i, \bm{x}_i)\}=&\bm{\alpha}_0(u_i)^T\bm{m}_i+\bm{\alpha}_1(u_i)^T\bm{x}_i=\sum_{j=1}^p \alpha_{0j}(u_i)m_{ij}+\sum_{l=1}^q\alpha_{1l}(u_i)x_{il},\\
    \bm{m}_i&=\bm{\Gamma}(u_i)^T\bm{x}_i+\bm{\epsilon}_i=\sum_{l=1}^q \bm{\Gamma}_l(u_i)x_{il}+\bm{\epsilon}_i.
\end{align*}

To explore potential effects of the varying covariate $u$, we apply the regression spline method to approximate coefficient functions. Using cubic B-splines, we approximate $\{\alpha_{0j}(u_i), j=1,\cdots, p\}$, $\{\alpha_{1l}(u_i), l=1,\cdots, q\}$ and \{$\Gamma_{lj}(u_i),l=1,\cdots,q, j=1,\cdots, p\}$ as follows:
\begin{eqnarray*}
&&\alpha_{0j}(u_i)\approx \sum_{k=1}^{K_0}a_{0jk}b_{0jk}(u_i),\quad
\alpha_{1l}(u_i)\approx  \sum_{k=1}^{K_1}a_{1lk}b_{1lk}(u_i),\\
&&\Gamma_{lj}(u_i)\approx \sum_{k=1}^{K_m}c_{ljk}B_{ljk}(u_i),
\end{eqnarray*}
where for any $j=1,\cdots, p$, $l=1,\cdots, q,$ $\{b_{0jk}(\cdot),k=1\cdots,K_0\}$, $\{b_{1lk}(\cdot),k=1\cdots,K_1\}$ and $\{B_{ljk}(\cdot), k=1,\cdots, K_m\}$ are sets of B-spline bases in linear spaces $\mathbb{G}_{\alpha_{0j}}$, $\mathbb{G}_{\alpha_{1 l}}$ and $\mathbb{G}_{\Gamma_{lj}}$ of spline functions on $\mathcal{U}$, respectively. Degrees and knots are fixed, which may be different across $k$. $K_0$, $K_1$ and $K_m$ are the numbers of basis functions used for $\alpha_{0j}(\cdot)$, $\alpha_{1l}(\cdot)$ and $\Gamma_{lj}(\cdot)$, $l=1,\cdots, q,$ $j=1,\cdots,p,$ respectively. $K_0$, $K_1$ and $K_m$ are finite numbers chosen to be large enough for accurate approximation. 

Since the $l$th element of $\bm{\beta}(u)$ is only determined by the $l$th row of $\bm{\Gamma}(u)$, we consider the basis expansion of $\bm{\Gamma}(u)$ row by row. Let $B_{lk}(u)\equiv B_{ljk}(u)$, for all $j=1,\cdots, p$, for simplicity. We express $\{\bm{\Gamma}_l(u_i), l=1,\cdots q\}$ as
$$\bm{\Gamma}_l(u_i)\approx\sum_{k=1}^{K_m}\bm{c}_{lk}B_{lk}(u_i)$$
where $\bm{\Gamma}_l(u_i)\in\mathbb{R}^{p\times 1}$ is the $l$th row of $\bm{\Gamma}(u_i)$, $\bm{\Gamma}_l(u_i)=(\Gamma_{l1},\cdots, \Gamma_{lp})^T$,
$\bm{c}_{lk}=(c_{l1k},\cdots, c_{lpk})^T\in\mathbb{R}^{p\times 1}$, and for any $l=1,\cdots, q,$ $\{B_{lk}(\cdot), k=1,\cdots, K_m\}$ is a basis of spline functions on $\mathcal{U}$ with a fixed degree and knots for a linear space $\mathbb{G}_{\Gamma_l},$ which is equivalent to $\mathbb{G}_{\Gamma_{lj}}$.
Then
\begin{align}
\label{GVCMMappox_y}g\{\mathrm{E}(y_i\mid u_i,\bm{m}_i, \bm{x}_i)\}&\approx \sum_{j=1}^p \sum_{k=1}^{K_0}a_{0jk}b_{0jk}(u_i)m_{ij}+\sum_{l=1}^q \sum_{k=1}^{K_1}a_{1lk}b_{1lk}(u_i)x_{il}\nonumber\\
&\overset{\triangle}=\sum_{j=1}^p \bm{a}_{0j}^T\bm{b}_{0j}(u_i)m_{ij}+\sum_{l=1}^q \bm{a}_{1l}^T\bm{b}_{1l}(u_i)x_{il},
\end{align}
and
\begin{align}
\label{GVCMMapprox_m} \bm{m}_i\approx \sum_{l=1}^q\sum_{k=1}^{K_m}\bm{c}_{lk}B_{lk}(u_i)x_{il}+\bm{\epsilon}_i 
\overset{\triangle}{=}\sum_{l=1}^q\bm{c}_l^T\bm{B}_l(u_i)x_{il}+\bm{\epsilon}_i.
\end{align}

Let $\bm{\alpha}_0^*=(\bm{a}_{01}^T,\cdots,\bm{a}_{0p}^T)^T\in\mathbb{R}^{pK_0\times1}$, $\bm{m}_i^*=(m_{i1}\bm{b}_{01}(u_i)^T,\cdots,m_{ip}\bm{b}_{0p}(u_i)^T)^T\in\mathbb{R}^{pK_0\times 1}$, $\bm{\alpha}_1^*=(\bm{a}_{11}^T,\cdots,\bm{a}_{1q}^T)^T\in\mathbb{R}^{qK_1\times1}$, $\bm{x}_i^*=(x_{i1}\bm{b}_{11}(u_i)^T,\cdots,x_{iq}\bm{b}_{1q}(u_i)^T)^T\in\mathbb{R}^{qK_1\times 1}$, $\bm{C}=(\bm{c}_1^T,\cdots,\bm{c}_q^T)^T\in\mathbb{R}^{qK_m\times p}$ and $\bm{x}_i^m=(x_{i1}\bm{B}_1(u_i)^T,\cdots, x_{iq}\bm{B}_q(u_i)^T)^T\in\mathbb{R}^{qK_m\times 1}$. Then equations (\ref{GVCMMappox_y}) and (\ref{GVCMMapprox_m}) become
$$g\{\mathrm{E}(y_i\mid u_i,\bm{m}_i, \bm{x}_i)\}\approx\bm{\alpha}_0^{*T}\bm{m}_i^*+\bm{\alpha}_1^{*T}\bm{x}_i^*,\ \text{and}
\ \bm{m}_i\approx\bm{C}^T\bm{x}_i^m+\bm{\epsilon}_i.$$

By Newton-Raphson algorithm or Fisher Scoring algorithm, we obtain the estimates of $\bm{\alpha}_0^*$ and $\bm{\alpha}_1^*$. The estimate of $\bm{C}$ can be obtained by minimizing $\ell(\bm{C})=\left\Vert\bm{M}-\bm{X}^m\bm{C}\right\Vert^2$ with respect to $\bm{C}$. That is,
$$\hat{\bm{C}}=(\bm{X}^{mT}\bm{X}^m)^{-1}\bm{X}^{mT}\bm{M},$$
where $\bm{M}=(\bm{m}_1,\cdots, \bm{m}_n)^T$, $\bm{X}^m=(\bm{x}_1^{m},\cdots \bm{x}_n^{m})^T$.

To obtain the estimates of $\bm{\alpha}_0(u)$, $\bm{\alpha}_1(u)$ and $\bm{\Gamma}(u)$, for ease of presentation, set $\bm{b}_{0j}(u_i)\equiv \bm{b}_{0}(u_i)$, for all $j=1,\cdots,p$, $\bm{b}_{1l}(u_i)\equiv \bm{b}_{1}(u_i)$, for all $l=1,\cdots,q$, and $\bm{B}_l(u_i)=\bm{B}(u_i)$, for all $l=1,\cdots q$.
Thus,
\begin{eqnarray*}
&&\hat{\bm{\alpha}}_0(u_i)=(\bm{b}_0(u_i)^T\hat{\bm{a}}_{01},\cdots \bm{b}_0(u_i)^T\hat{\bm{a}}_{0p})^T=\{\bm{I}_p\otimes \bm{b}_0^T(u_i)\}\hat{\bm{\alpha}}_0^*,\\
&&\hat{\bm{\alpha}}_1(u_i)=(\bm{b}_1(u_i)^T\hat{\bm{a}}_{11},\cdots \bm{b}_1(u_i)^T\hat{\bm{a}}_{1p})^T=\{\bm{I}_q\otimes \bm{b}_1^T(u_i)\}\hat{\bm{\alpha}}_1^*,\\
&&\hat{\bm{\Gamma}}(u_i)=(\hat{\bm{c}}_1^T\bm{B}(u_i),\cdots, \hat{\bm{c}}_q^T\bm{B}(u_i))^T=\{\bm{I}_q\otimes \bm{B}(u_i)^T\}\hat{\bm{C}}.
\end{eqnarray*}
The estimates $\hat{\bm{\alpha}}_0(u)$ and $\hat{\bm{\Gamma}}(u)$ can be obtained at any point $u$ based on basis functions. Naturally, we estimate the direct $\bm{\alpha}_1(u)$ and indirect effects $\bm{\beta}(u)$ by $\hat{\bm{\alpha}}_1(u)$ and $\hat{\bm{\beta}}(u)=\hat{\bm{\Gamma}}(u)\hat{\bm{\alpha}}_0(u)$. \blue{In practice, one needs to determine $K_0$, $K_1$ and $K_m$, which control the model complexity of $\hat{\bm{\alpha}}_0(u)$, $\hat{\bm{\alpha}}_1(u)$ and $\hat{\bm{\Gamma}}(u)$, respectively. Cross-validation may be used to select $K_0$, $K_1$ and $K_m$. For ease of computation, we would suggest setting $K_0$, $K_1$ and $K_m$ to be the same in implementation.}

\subsection{Asymptotic Theory}
In this section, we derive the asymptotic properties of the estimated direct effect $\hat{\bm{\alpha}}_1(u)$ and the estimated indirect effect $\hat{\bm{\beta}}(u)$. First note that
$$\bm{\alpha}_0(u_i)\approx(\bm{b}_0(u_i)^T\bm{a}_{01},\cdots \bm{b}_0(u_i)^T\bm{a}_{0p})^T=\{\bm{I}_p\otimes \bm{b}_0^T(u_i)\}\bm{\alpha}_0^*,$$
$$\bm{\alpha}_1(u_i)\approx(\bm{b}_1(u_i)^T\bm{a}_{11},\cdots \bm{b}_1(u_i)^T\bm{a}_{1p})^T=\{\bm{I}_q\otimes \bm{b}_1^T(u_i)\}\bm{\alpha}_1^*,$$
$$\bm{\Gamma}(u_i)\approx(\bm{c}_1^T\bm{B}(u_i),\cdots, \bm{c}_q^T\bm{B}(u_i))^T=\{\bm{I}_q\otimes \bm{B}(u_i)^T\}\bm{C}.$$

The following technical conditions along with some notations are imposed.

\begin{itemize}
    \item [\textbf{C1.}] The sample $\{u_{i}, i=1,\cdots, n\}$ is independently distributed with distribution $F_U$ on a bounded support $\mathcal{U}$ and Lebesgue density $f_U(u)$ which is bounded away from 0 and infinity uniformly over $u\in \mathcal{U}.$

    \item [\textbf{C2.}]
    The functions $\{\alpha_{0j}(u)\}_{j=1}^p$ and $\{\alpha_{1l}(u)\}_{l=1}^q$ belong to a class of functions $\mathcal{F}_1$, whose $r_1$th derivatives $\alpha_{0j}^{(r_1)}$ and $\alpha_{1l}^{(r_1)}$ exist and are Lipschitz of order $\eta_1$,
    $$\mathcal{F}_1=\{\alpha_j(\cdot):\mid\alpha_j^{(r_1)}(s)-\alpha_j^{(r_1)}(u)\mid \leq K_1\mid s-t\mid ^{\eta_1} \text{ for } s,t\in \mathcal{U}\},$$
    for some positive constant $K_1$, where $r_1$ is a nonnegative integer and $\eta_1\in(0,1]$ such that $\nu_1=r_1+\eta_1>0.5.$\\

    Similarly, the functions $\{\Gamma_{lj}(u)\}$, for any $j=1,\cdots, p$, $l=1,\cdots q$, belong to a class of functions $\mathcal{F}_2$, whose $r_2$th derivatives $\Gamma_{lj}^{(r_2)}$ exists and are Lipschitz of order $\eta_2$,
    $$\mathcal{F}_2=\{\Gamma_{lj}(\cdot): \mid \Gamma_{lj}^{(r_2)}(s)-\Gamma_{lj}^{(r_2)}(u)\mid \leq K_2\mid s-t\mid ^{\eta_2} \text{ for } s,t\in \mathcal{U}\},$$
    for some positive constant $K_2$, where $r_2$ is a nonnegative integer and $\eta_2\in(0,1]$ such that $\nu_2=r_2+\eta_2>0.5.$

    \item [\textbf{C3.}] The eigenvalues $\lambda_1\leq \cdots \leq \lambda_q$ of $\bm{\Sigma}_{xx}=\mathrm{E}(\bm{x}_i\bm{x}_i^T)$ and eigenvalues $\lambda'_1\leq \cdots \leq \lambda'_p$ of $\bm{\Sigma}_{mm}=\mathrm{E}(\bm{m}_i\bm{m}_i^T)$
    are bounded away from 0 and infinity. That is, $N_1\leq \lambda_1\leq \cdots \leq \lambda_q\leq N_2$ and $N'_1\leq \lambda'_1\leq \cdots \leq \lambda'_p\leq N'_2$ for some positive constants $N_1$, $N_2$, $N'_1$, and $N'_2$.

    \item [\textbf{C4.}] All eigenvalues of $\Sigma_{\bm{\epsilon}}$ are bounded away from 0 and infinity.

    \item [\textbf{C5.}] The second-order derivatives of coefficients $\bm{\Gamma}(u)$, $\bm{\alpha}_0(u)$ and $\bm{\alpha}_1(u)$ are assumed to be continuous over $\mathcal{U}$. Thus, they and their second-order derivatives are bounded. Denote $M_\Gamma$, $M_0$ and $M_1$ to be the bounds of $\bm{\Gamma}(u)$, $\bm{\alpha}_0(u)$ and $\bm{\alpha}_1(u)$ over $u\in\mathcal{U}$, respectively.

    \item[\textbf{C6.}] $\limsup_n\bigg(\frac{\max_\omega K_\omega}{\min_\omega K_\omega}\bigg)<\infty$, $\omega = 0, 1, m$.

    \item[\textbf{C7.}] There exist a positive constant $N_3$ such that $\mid x_{il}\mid \leq N_3$, for all $i=1,\cdots, n$, $l=1,\cdots q$.

    \item [\textbf{C8.}] There exist constants $C_1>0$, such that for sufficiently large $n$,
    $$C_1\leq \lambda_{min}\{-n^{-1}L^{''}(\bm{\alpha}_0, \bm{\alpha}_1)\}$$
    for $(\bm{\alpha}_0,\bm{\alpha}_1)\in \cup (\bm{\alpha}^*_0,\bm{\alpha}^*_1)$, where $L(\bm{\alpha}^*_0,\bm{\alpha}^*_1)$ is the log-likelihood function based on equation (\ref{GVCMMmodel_y}), and $\lambda_{min}(\cdot)$ denotes the smallest eigenvalue of a matrix.

\end{itemize}

Conditions \textbf{C1} and \textbf{C2} ensure the property of
B-spline approximation. In particular, \textbf{C1} guarantees that
the observations $u_i$, $i=1,\cdots, n$, are randomly scattered.
Condition \textbf{C2}, adapted from \cite{yang2017feature},
guarantees the smoothness of the coefficient functions. Condition
\textbf{C3} and \textbf{C4} ensure non-singularity of the
covariance matrices. Condition \textbf{C5} is required to derive
the convergence rate of $\hat{\bm{\Gamma}}(u)$,
$\hat{\bm{\alpha}}_0(u)$ and $\hat{\bm{\beta}}(u)$. Conditions
\textbf{C6-C8}, adapted from \cite{huang2004polynomial} and
\cite{yang2017feature}, are technical conditions to facilitate
theoretical understanding of the proposed estimation procedure.
Condition \textbf{C8} is needed to control the second derivative
of the log-likelihood function in the derivation of convergence
rate of $\hat{\bm{\beta}}(\cdot)$. These are some mild conditions
that can be satisfied in many practical scenarios. In practice, the approximation error of
spline expansion becomes very small and therefore can be negligible when the number of knots
$K_0$, $K_1$ and $K_m$ are large enough. Thus, it is assumed throughout this paper that $\alpha_{0j}(u)$,
$\alpha_{1l}(u)$ and $\Gamma_{lj}(u)$ belong to the linear spaces
$\mathbb{G}_{\alpha_{0j}}$, $\mathbb{G}_{\alpha_{1 l}}$ and
$\mathbb{G}_{\Gamma_{lj}}$, respectively, for technical
simplicity. Thus, the approximation error due to spline approximation may be ignored.

 Let $\left\Vert a(\cdot)\right \Vert_{L_2}$ denote the $L_2$ norm of a square integrable function $a(u)$ on $\mathcal{U}$, i.e. $\left\Vert a(\cdot)\right\Vert_{L_2}^2=\int_{u\in\mathcal{U}}\mid a(u)\mid ^2 f_U(u) du$, and
 \begin{align*}
    \bm{\Sigma}_{m^*m^*}&=\mathrm{E}\{\bm{m}_i^*\bm{m}_i^{*T}b^{''}(\bm{\alpha}_0^{*T}\bm{m}_i^*+\bm{\alpha}_1^{*T}\bm{x}_i^*)\},\\
    \bm{\Sigma}_{m^*x^*}&=\mathrm{E}\{\bm{m}_i^*\bm{x}_i^{*T}b^{''}(\bm{\alpha}_0^{*T}\bm{m}_i^*+\bm{\alpha}_1^{*T}\bm{x}_i^*)\},\\
    \bm{\Sigma}_{x^*m^*}&=\mathrm{E}\{\bm{x}_i^*\bm{m}_i^{*T}b^{''}(\bm{\alpha}_0^{*T}\bm{m}_i^*+\bm{\alpha}_1^{*T}\bm{x}_i^*)\},\\
    \bm{\Sigma}_{x^*x^*}&=\mathrm{E}\{\bm{x}_i^*\bm{x}_i^{*T}b^{''}(\bm{\alpha}_0^{*T}\bm{m}_i^*+\bm{\alpha}_1^{*T}\bm{x}_i^*)\}.
\end{align*}

The asymptotic distribution of the estimated direct effect $\hat{\bm{\alpha}}_1(u)$ can be summarized as follows.
\begin{theorem}\label{GVCMMalpha1_u}
Under Conditions \textbf{C1-C6}, the estimated direct effect $\hat{\bm{\alpha}}_1(u)$ satisfies that
\begin{align*}
[\Var\{\hat{\bm{\alpha}}_1(u)\}]^{-1/2}\{\hat{\bm{\alpha}}_1(u)-\bm{\alpha}_1(u)\}\overset{D}{\rightarrow}N(\bm{0},\bm{I}),
\end{align*}
where $\Var\{\hat{\bm{\alpha}}_1(u)\} = \frac{1}{n}\phi\{\bm{I}_q\otimes \bm{b}_1(u)^T\}\bm{\Sigma}_{x^*x^*\mid m^*}^{-1}\{\bm{I}_q\otimes \bm{b}_1(u)\}$, $\bm{\Sigma}_{x^*x^*\mid m^*}=\bm{\Sigma}_{x^*x^*}-\bm{\Sigma}_{x^*m^*}\bm{\Sigma}_{m^*m^*}^{-1}\bm{\Sigma}_{m^*x^*}$, and $\otimes$ denotes the Kronecker product.
\end{theorem}

The asymptotic normality sets up basis for subsequentially statistical inference towards $\bm{\alpha}_1(u)$. To obtain the asymptotic distribution of $\hat{\bm{\beta}}(u) = \hat{\bm{\Gamma}}(u)\hat{\bm{\alpha}}_0(u)$, we consider the following two lemmas at first.

\begin{lemma}\label{GVCMMlem_alpha0_u}
Under Conditions \textbf{C1-C6}, the estimate $\hat{\bm{\alpha}}_0(t)$ satisfies
\begin{align*}
[\Var\{\hat{\bm{\alpha}}_0(u)\}]^{-1/2}\{\hat{\bm{\alpha}}_0(u)-\bm{\alpha}_0(u)\}\overset{D}{\rightarrow}N(\bm{0},\bm{I}),
\end{align*}
where $\Var\{\hat{\bm{\alpha}}_0(u)\} = \frac{1}{n}\phi\{\bm{I}_p\otimes \bm{b}_0(u)^T\}\bm{\Sigma}_{m^*m^*\mid x^*}^{-1}\{\bm{I}_p\otimes \bm{b}_0(u)\}$, $\bm{\Sigma}_{m^*m^*\mid x^*}=\bm{\Sigma}_{m^*m^*}-\bm{\Sigma}_{m^*x^*}\bm{\Sigma}_{x^*x^*}^{-1}\bm{\Sigma}_{x^*m^*}$, and $\otimes$ denotes the Kronecker product.
\end{lemma}

Lemma 1 establishes the asymptotic properties of $\hat{\bm{\alpha}}_0(u)$, and Lemma 2 establishes the asymptotic properties of $\hat{\bm{\Gamma}}(u)$.

\begin{lemma}\label{GVCMMlem_Gamma_u}
Under Conditions \textbf{C1-C6}, the estimate of $\hat{\bm{\Gamma}}(u)$ is
$$\hat{\bm{\Gamma}}(u)=\{\bm{I}_q\otimes \bm{B}(u)^T\}(\bm{X}^{mT}\bm{X}^m)^{-1}\bm{X}^{mT}\bm{M}$$
and then
$$[\Var\{\hat{\bm{\Gamma}}(u)\}]^{-1/2}\{\hat{\bm{\Gamma}}(u)-\bm{\Gamma}(u)\}\overset{D}{\rightarrow}N(\bm{0},\bm{I}),$$
where $\Var\{\hat{\bm{\Gamma}}(u)\}=\frac{1}{n}\big[\{\bm{I}_q\otimes \bm{B}(u)^T\}\bm{\Sigma}_{x^mx^m}^{-1}\{\bm{I}_q\otimes \bm{B}(u)\}\big]\otimes\bm{\Sigma}_\epsilon$, $\otimes$ denotes the Kronecker product, $\bm{M}=(\bm{m}_1,\cdots, \bm{m}_n)^T$, $\bm{X}^m=(\bm{x}_1^{m},\cdots,$ $\bm{x}_n^{m})^T$, and $\bm{\Sigma}_{x^mx^m}=\mathrm{E}{(\bm{x}_i^m\bm{x}_i^{mT})}$.
\end{lemma}

We next study the asymptotic properties of the indirect effect estimate $\hat{\bm{\beta}}$.

\begin{theorem}[Convergence Rate]\label{GVCMMrates} Under Conditions \textbf{C1-C8}, $\left\Vert\hat{\beta}_l-\beta_l\right\Vert^2_{L_2}=O_p(1/n).$
\end{theorem}

Note that $\hat{\bm{\alpha}}_0(u)$ and $\hat{\bm{\Gamma}}(u)$ are asymptotically unbiased estimate of $\bm{\alpha}_0(u)$ and $\bm{\Gamma}(u)$, respectively. Assuming the coefficient functions belong to linear spaces of B-splines, the bias term can be ignored. Therefore, $\hat{\beta}_l-\beta_l$ contributes mainly to the variance of estimation. The asymptotic distribution of $\hat{\bm{\beta}}(u)$ is given in the following theorem.

\begin{theorem}[Asymptotic Distribution]\label{GVCMMasy.} Under Conditions \textbf{C1-C6}, $$[\Var\{\hat{\bm{\beta}}(u)\}]^{-1/2}\{\hat{\bm{\beta}}(u)-\bm{\beta}(u)\}\overset{D}{\rightarrow}N\{\bm{0},\bm{I}\},$$
where $\Var\{\hat{\bm{\beta}}(u)\}=\frac{1}{n}\bm{\alpha}_0(u)^T\Sigep \bm{\alpha}_0(u) \{\bm{I}_q\otimes \bm{B}(u)^T\}\bm{\Sigma}_{x^mx^m}^{-1}\{\bm{I}_q\otimes \bm{B}(u)\}+\frac{\phi}{n}\bm{\Gamma}(u)\{\bm{I}_p\otimes \bm{b}_0(u)^T\}\bm{\Sigma}_{m*m*\mid x*}^{-1}\{\bm{I}_p\otimes \bm{b}_0(u)\}\bm{\Gamma}(u)^T$ and $\bm{\Sigma}_{m^*m^*\mid x^*}=\bm{\Sigma}_{m^*m^*}-\bm{\Sigma}_{m^*x^*}\bm{\Sigma}_{x^*x^*}^{-1}\bm{\Sigma}_{x^*m^*}$.

\end{theorem}

\section{Simulation Studies}\label{sec3}

In this section, we conduct simulation studies to examine the performance of the proposed procedures. Different sample sizes $n=300$ and $n=500$ are considered. The covariate $u-$ and $\bm{x}-$ variables are generated as follows. First, draw $(u_i^*,\bm{x}_{-1,i})$ from a $q=3$ dimensional normal distribution $N_3(\bm{0},\bm{\Sigma})$, where the $(k_1,k_2)$th element of $\bm{\Sigma}$ is set to be $\rho^{\mid k_1-k_2\mid }$ with $\rho=0.5$. Then set $u_i=\Phi(u_i^*)$ where $\Phi(\cdot)$ is the cumulative distribution function of the standard normal distribution. Therefore, $u_i$ is uniformly distributed on $[0,1]$ and correlated with $\bm{x}_{-1,i}.$ An additional column of 1's is added to $\bm{x}_{-1,i}$ for intercept, so that the dimension of $\bm{x}$ is equal to $q=3$.

One example is designed for a binary response and the other for a Poisson response. For all following settings, we denote $\bm{\alpha}_0(u)=(\alpha_{01}(u),\alpha_{02}(u))^T$, $\bm{\alpha}_1(u)=(\alpha_{11}(u),\alpha_{12}(u),\alpha_{13}(u))^T$, and
\begin{equation*}
\bm{\Gamma}(u)=\left(
 \begin{array}{cc}
    \Gamma_{11}(u) & \Gamma_{12}(u)\\
    \Gamma_{21}(u) & \Gamma_{22}(u)\\
    \Gamma_{31}(u) & \Gamma_{32}(u)\\
 \end{array}
\right).
\end{equation*}

$\emph{Example \uppercase\expandafter{\romannumeral1} (Binary response)}.$ Based on model (\ref{GVCMMmodel_m}), for $i=1,\cdots, n$, we generate mediators from $\bm{m}_i=\bm{\Gamma}^T(u_i)\bm{x}_i+\bm{\epsilon}_i$, where $\bm{\epsilon}\sim N(0,0.5\bm{I})$. Then, we generate binary responses with the probability of $y_i=1$ being $p(u_i,\bm{x}_i)$, $i=1,\cdots, n$, defined below:
$$\text{logit}\{p(u_i,\bm{x}_i)\}=\bm{\alpha}_0(u_i)^T\bm{m}_i+\bm{\alpha}_1(u_i)^T\bm{x}_i,$$
where $\text{logit}\{t\}=\log\{t/(1-t)\}$, the logit link in the logistic regression.

The coefficient functions are defined by
\begin{align*}
\alpha_{01}(u)=&~u(1-u), \quad \alpha_{02}(u)=\sin(\pi u) + 0.5,\\
\alpha_{11}(u)=&~0.5\Phi_1(u), \quad  \alpha_{12}(u)=0.5u^3, \quad  \alpha_{13}(u)=0.5\Phi_2(u),\\
\Gamma_{11}(u)=&~2\sin(\pi u), \quad \Gamma_{12}(u)=(1-u)^2, \\
\Gamma_{21}(u)=&~0.5(0.5-u)^2,  \quad \Gamma_{22}(u)=0.5\cos(\pi u -\frac{\pi}{2}), \\
\Gamma_{31}(u)=&~-\cos(\pi u-\frac{\pi}{2}),  \quad \Gamma_{32}(u)=3(u-0.5)^2,
\end{align*}
where $\Phi_1(\cdot)$ and $\Phi_2(\cdot)$ are the cumulative distribution functions of $Gamma(0.5,1)$ and $N(0,1)$, respectively.

$\emph{Example \uppercase\expandafter{\romannumeral2} (Poisson response)}.$ For $i=1,\cdots, n$, we generate mediators from $\bm{m}_i=\bm{\Gamma}^T(u_i)\bm{x}_i+\bm{\epsilon}_i$, where $\bm{\epsilon}\sim N(0,0.25\bm{I})$. Then, given $\{u_i,\bm{x}_i\}$, $i=1,\cdots, n$, we generate count responses $y_i$ from Poisson distributions with mean $\lambda\{u_i, \bm{x}_i\}$ defined below:
$$\log\{\lambda(u_i,\bm{x}_i)\}=\bm{\alpha}_0(u_i)^T\bm{m}_i+\bm{\alpha}_1(u_i)^T\bm{x}_i.$$

The coefficient functions are defined by
\begin{align*}
\alpha_{01}(u)=&~-0.05e^{2u-1}-0.5, \quad \alpha_{02}(u)=-0.05u(0.25-u)-1,\\
\alpha_{11}(u)=&~0.02\sin^2(2\pi u), \quad  \alpha_{12}(u)=0.02\cos^2(\pi u), \\ \alpha_{13}(u)=&0.05 u^2 ,\\
\Gamma_{11}(u)=&~0.1u^2, \quad \Gamma_{12(u)}=0.2(e^u+e^{-u})-0.5, \\
\Gamma_{21}(u)=&~0.1\cos(\pi u - \frac{\pi}{2}),  \quad \Gamma_{22}(u)=0.3\sin^2(\pi u), \\
\Gamma_{31}(u)=&~-0.05u^4+0.02,  \quad \Gamma_{32}(u)=0.5\sin(\pi u) + 0.02.
\end{align*}

We employ B splines for approximation, and set the degree of splines to be three, indicating the most commonly used cubic spline.  The number of interior knots is set to be three and the intercept is considered, giving rise to a total of seven degrees of freedom  for each varying coefficient. Since the variable $u$ is almost uniform over $[0,1]$, we employ equally spaced knots. For simplicity, we set the basis splines for approximations of $\alp(u)$, $\alpp(u)$ and $\Gam(u)$ the same in the study. That is, $\bm{b}_0(u_i)=\bm{b}_1(u_i)=\bm{B}(u_i)$, $i=1,\cdots,n$. \blue{Different B-spline basis functions can be employed for different coefficient functions based on prior knowledge and expertise. In practice, this is a subject and empirically-driven task. For example, knots can be placed at locations with changing curvatures of coefficients, and the degrees of B-spline can be determined based on smoothness of coefficients.}

To interpret the final results more clearly, instead of at randomly distributed points, $\hat{\bm{\alpha}}_0(u)$ and $\hat{\Gam}(u)$ are estimated at a set of grid points $\{u_k,k=1,\cdots, n_{grid}\}$ between 0 and 1, with an increment of 0.002 and a total of $n_{grid}=500$ points. Once we have the estimated $\bm{\alpha}_0^*$ and $\bm{C}$, the estimation results and corresponding asymptotic variance at each grid point can be derived.

\begin{figure}[h]
\centering
\includegraphics[width=10cm,height=9cm]{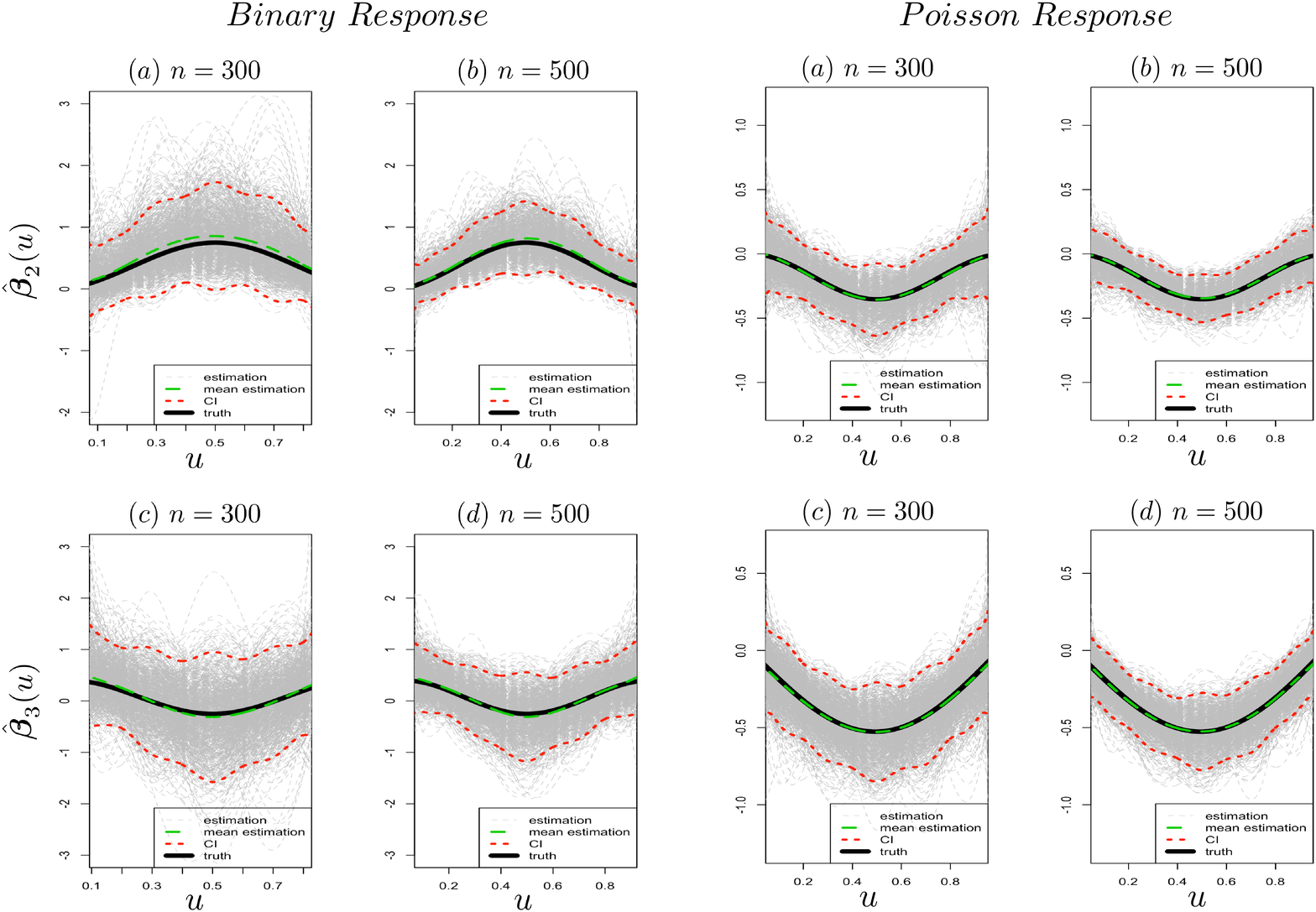}
\caption[The Mean of Estimated Coefficient Functions]{The Mean of Estimated Coefficient Functions. The red and green curves are the  estimated values and true values, respectively. The dashed curves are the estimated function plus/minus 1.96 times standard errors in 1000 repetitions.}
\label{SimCombo}
\end{figure}

Figure \ref{SimCombo} depicts the estimated coefficient functions at all grid points for binary response and Poisson response, respectively. The means of the estimated coefficient functions are plotted in green dashed curves, while the true values are in black. Dashed red curves denote the 95\% pointwise confidence bands based on 1000 simulations. As can be seen from Figure \ref{SimCombo}, the proposed estimation procedure gives a reasonably good approximation to the true coefficient functions. The true curves falls in the corresponding 95\% confidence bands. As $n$ increases, the confidence interval becomes narrower.

The accuracy of the proposed standard error estimation procedure and coverage probability of true curves are examined in the following tables. Table \ref{GVCMMtable_setting1} summarizes the simulation results for binary response with the sample sizes $n=300$ and $n=500$, and Table \ref{GVCMMtable_setting2} summarizes the simulation results for poisson response. \emph{Mean(sd)} in the first column is the mean and standard deviation of $\hat{\beta}_2(u)$ and $\hat{\beta}_3(u)$, excluding the element of $\hat{\bm{\beta}}(t)$ with respect to the intercept column in $\bm{X}$, at different points of $u$ in 1000 simulations. The standard deviations of estimates in 1000 repetitions can be considered as the true deviations of estimates. \emph{SE(sd)} in the second column denotes the average and standard deviation of 1000 estimated standard errors, calculated based on Theorem 3. As can be seen, the average standard errors in the second column is very close to the standard deviation in the first column. This implies that the standard error estimated based on Theorem 3 performs well. Actually, all differences between average standard errors and standard deviations are within one standard deviation, for both binary and poisson responses. \emph{Coverage probability (CP)} is the probability of true values being covered by corresponding confidence intervals at the significance level of 95\%. That is, estimated functions plus/minus 1.96 times estimated standard errors. The coverage probabilities are all around 0.95. \blue{The Monte Carlo errors for 1000 simulations is 1.35\% for confidence level 95\%. As can be seen in Tables 1 and 2, most CP values lie in 95\% $\pm$ 1.35\%. This implies that Theorem 3 are valid.}

\begin{table}[htbp]
\caption{Simulation Results of Example \uppercase\expandafter{\romannumeral1} for Binary Response}
\label{GVCMMtable_setting1}
\centering
\scalebox{1}{
\begin{tabular}{ccccccccccc}
    \toprule
    &&  \multicolumn{3}{c}{$\hat{\beta}_2(u)$}  &&\multicolumn{3}{c}{$\hat{\beta}_3(u)$}\\
    \cmidrule{3-5} \cmidrule{7-9}
    $n$&$u$ & mean(sd) & SE(sd) &CP && mean(sd) &SE(sd) &CP \\

    \midrule
    &1/4 &  .521(.336) & .317(.112) & .940 && .126(.468) & .436(.105)&.966\\
    300&1/2 &  .856(.447) & .408(.122) &.943 && -.311(.645)&.587(.115)&.945\\
    &3/4 &  .517(.365) & .325(.116) &.936&& .106(.496)&.458(.110)&.957\\
    &&&&&&&&\\
    &1/4 &  .469(.222) & .213(.055) & .943 && .118(.299) &.303(.049)&.958 \\
    500&1/2 &  .819(.306)  &  .294(.066)  & .952 &&  -.308 (.441)  & .420 (.060)  & .950\\
    &3/4 &    .479(.237)   &  .221 (.060)   & .933 &&  .102 (.331)  & .315 (.054)  & .956\\

    \bottomrule
\end{tabular}}
\end{table}

\begin{table}[htbp]
\caption{Simulation Results of Example \uppercase\expandafter{\romannumeral2} for Poisson Response}
\label{GVCMMtable_setting2}
\centering
\scalebox{1}{
\begin{tabular}{ccccccccccc}
    \toprule
    &&  \multicolumn{3}{c}{$\hat{\beta}_2(u)$}  &&\multicolumn{3}{c}{$\hat{\beta}_3(u)$}\\
    \cmidrule{3-5} \cmidrule{7-9}
    $n$&$u$ & mean(sd) & SE(sd) &CP && mean(sd) &SE(sd) &CP \\

    \midrule
    &1/4 &  -.195(.117) & .115(.025) & .946 && -.389(.134) & .133(.024)&.949\\
    300&1/2 &  -.356(.144) & .134(.028) &.935 && -.529(.164)&.162(.029)&.950\\
    &3/4 &  -.194(.118) & .121(.028) &.955&& -.373(.147)&.146(.030)&.949\\
    &&&&&&&&\\
    &1/4 &  -.185(.086) & .085(.014) & .949 && -.387(.099) &.099(.013)&.944 \\
    500&1/2 &  -.345(.095)  &  .098(.015)  & .945 &&  -.525 (.128)  & .120 (.016)  & .940\\
    &3/4 &    -.190(.096)   &  .089(.017)   & .940 &&  -.373 (.112)  & .108 (.017)  & .941\\

    \bottomrule
\end{tabular}}
\end{table}

\section{Real Data Example}\label{sec4}
We now illustrate the proposed methodology via an empirical analysis about the behavioral economics data. 300 small business owners were recruited to participate in a study by completing a survey online. Since only 262 individuals provided responses to all of the primary variables, the study is based on 262 responses (37 percent women). Entrepreneurial withdrawal intentions is the response variable ($y$), with two values - high and low. Economic stress ranging from 1 to 7 is considered to be the explanatory variable ($x$), where a larger number indicates increased economic stress. We treat the depressed affect as the mediator ($m$), which ranges from 1 to 5. Higher numbers indicate increased depressed affect.

\begin{figure}[h]
\centering
\includegraphics[width=8cm,height=5cm]{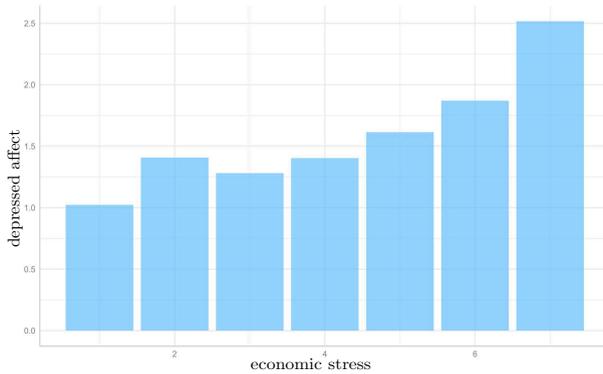}
\caption[Scatter Plot of Economic Stress and Depressed Affect]{Scatter Plot of Economic Stress and Depressed Affect}
\label{EDA_XM}
\end{figure}

In stressful economic times, when struggling to maintain and grow their business, entrepreneurs might develop feelings such as depression, helplessness and hopelessness as a result of job-related insecurity \citep{baum1986unemployment, dekker1995effects, hellgren2003does}. For instance, using our dataset, Figure \ref{EDA_XM} demonstrates that as economic stress increases, business owners are more likely to express higher depression levels. Economic stress and depressed affect are positively correlated. Thus, the explanatory effect of economic stress on depressed affect is of great interest.

\begin{figure}[h]
\centering
\includegraphics[width=8cm,height=5cm]{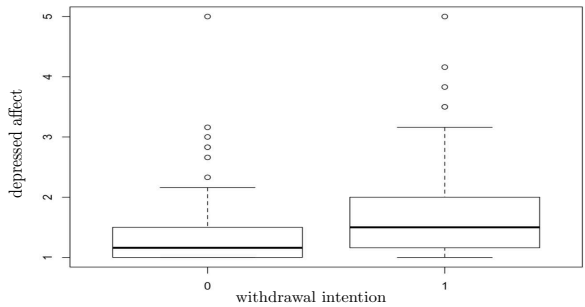}
\caption[Box Plot of Depressed Affect for Different Withdrawal Intention Groups]{Box Plot of Depressed Affect for Different Withdrawal Intention Groups}
\label{EDA_MY}
\end{figure}

Meanwhile, consistent with the learned helplessness theory, when the negative feelings are beyond the control, individuals are likely to withdraw from future entrepreneurship opportunities. This may due to that entrepreneurs' affect is a critical predictor of their engagement and withdrawal intentions. Figure \ref{EDA_MY} illustrates that higher intentions to withdraw are associated with higher level of depression. When entrepreneurs are unable to stay motivated during tough times, psychological withdrawal is inevitable.

Thus it is worthwhile exploring the mechanism of such chaining relationship. \cite{pollack2012moderating} studied the indirect effects of economic stress on entrepreneurs' intentions to withdraw through depressed affect based on a linear mediation model. It further shows that contact with business-related social ties plays as a moderator, buffering the impact of economic stress on depressed affect and in turn alleviating the potentially deleterious effects of economic stress. The relation between economic stress and withdrawal intentions was stronger among individuals with fewer social ties and weaker among those who reported more job-related contacts.

In addition to social ties, the moderating role of self-efficacy is also of interest. Self-efficacy is a person's belief in their ability to successfully reach their goal as a result of their own actions \citep{bandura1997self}. It can influence how people cope with stressors and make decisions. According to \cite{sequeira2007influence}, high self-efficacy is essential to overcoming substance abuse, achieving high academic performance and forming entrepreneurial intentions. Strong self-efficacy can provide critical confidence and support, alleviate feelings of depressed affect, and in turn decrease entrepreneurs' intentions to withdraw. Motivated by this hypothesis, we examine the moderating role of self-efficacy in the aforementioned mediation relationship. Instead of interaction terms, we explore the direct and indirect effects of economic stress as flexible functions of self-efficacy in this paper, so that the effects can be captured thoroughly and generalized well. That is, $u$ is defined to be self-efficacy, which is a continuous variable assessed based on \cite{chen1998does} 15-item scale.

In this analysis, we also include an intercept term ($z_1$) and other covariates such as age ($z_2$), tenure ($z_3$, how long an individual has worked in the company), social ties ($z_4$, the number of business-related contacts an entrepreneur has per day), social competence ($z_5$) and gender ($z_6$, male = 1, female = 0). Social ties is log-transformed to get rid of the influence of skewness.  Therefore, the model can be built up as
\begin{align}
\label{real_mod1}g\{\mathrm{E}(y_i\mid u,m, \bm{x})\}&=\alpha_0(u) m +  \alpha_1(u) x + \sum_{l=1}^6\alpha_{l+1}(u) z_{l}\\
\label{real_mod2}m&=\gamma_1(u) x + \sum_{l=1}^6\gamma_{l+1}(u) z_l + \epsilon.
\end{align}
For flexibility, we allow all coefficients to be smooth functions at first. Then to increase precision, we test significance of the varying pattern using procedures suggested by \cite{cai2000efficient}. Given varying $\gamma_1(\cdot)$, $\gamma_6(\cdot)$ and $\gamma_7(\cdot)$, the p-value for the null hypothesis that $\gamma_k(\cdot)$s are constant for $k=2,3,4,5$, is 0.712, illustrating constant effects. Thus, we set $\gamma_2$, $\gamma_3$, $\gamma_4$ and $\gamma_5$ to be constant. Furthermore, the p-value for the null hypothesis that $\alpha_q(\cdot)$s are constant is close to 1, $q=1,\cdots, 7$, and the p-value for the null hypothesis that $\alpha_0(\cdot)$ is constant is equal to 0.105, indicating a favor of constant effects. Thus, we consider the following model for subsequent analysis.
\begin{align}
\label{real_mod1}g\{\mathrm{E}(y_i\mid u,m, \bm{x})\}&=\alpha_0 m +  \alpha_1 x + \sum_{l=1}^6\alpha_{l+1} z_{l}\\
\label{real_mod2}m&=\gamma_1(u) x + \sum_{l=1}^4\gamma_{l+1} z_l + \sum_{l=5}^6\gamma_{l+1}(u) z_l + \epsilon.
\end{align}

\begin{figure}[h]
\centering
\includegraphics[width=10cm,height=6cm]{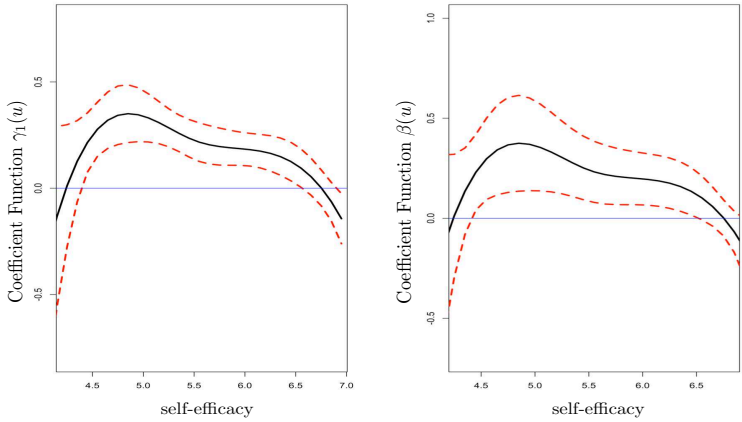}
\caption[The Estimated Coefficient Functions for $\gamma_1(u)$ and the Indirect Effect $\beta(u)$ for Economic Stress]{The Estimated Coefficient Functions for $\gamma_1(u)$ and the Indirect Effect $\beta(u)$ for Economic Stress. The dashed curves are the estimated functions plus/minus 1.96 times the estimated standard errors.}
\label{fig_real1}
\end{figure}

Cubic splines are employed in the analysis. The number of knots is chosen to be 1 by 10-fold cross validation. The left panel of Figure \ref{fig_real1} shows the effect of economic stress on depressed affect $\gamma_1(\cdot)$. The red dashed curves denote the 95\% point-wise confidence band. As can be seen, the effect is positive and changes over self-efficacy. Individuals develop depressed affect under economic stress, but self-efficacy serves to buffer the impact of economic stress on depressed affect when it is high. The small value of $\gamma_1$ when self-efficacy is low might be explained by the fact that in the start-up process, some participants who are less confident about themselves are still energetic and optimistic about potential deleterious effects of economic stress. The confidence interval on the left side is wide since a few participants reported extremely low scores. The effect of depressed affect on withdrawal intentions $\alpha_0$ is estimated to be constant at 1.068, significant at the 0.05 level. Thus, the indirect effect of economic stress on entrepreneurs' intentions to withdraw $\beta(\cdot)=\gamma_1(\cdot)\alpha_0$ shares the same pattern as $\gamma_1(\cdot)$, illustrated in the right panel of Figure \ref{fig_real1}. The direct effect of economic stress on withdrawal intentions is not significant.

\begin{figure}[h]
\centering
\includegraphics[width=6cm,height=6cm]{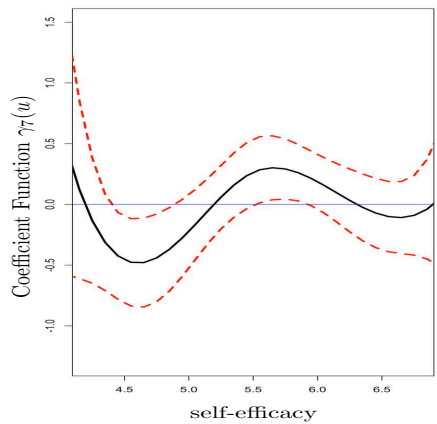}
\caption[The Estimated Coefficient Functions for $\gamma_7(u)$ for Gender]{The Estimated Coefficient Functions for $\gamma_7(u)$ for Gender. The dashed curves are the estimated functions plus/minus 1.96 times the estimated standard errors.}
\label{fig_real2}
\end{figure}

The effect of gender on depressed affect is also interesting. As illustrated in Figure \ref{fig_real2}, when self-efficacy is low, compared with females, males are less likely to have depressed effect. This might be explained by stronger ability to withstand pressure in the early stage of business for males. When self-efficacy score is high, the difference between males and females is not obvious.

In comparison, a simple mediation analysis is conducted based on linear models, with no effect modifier taken into account. The effect of economic stress on depressed affect $\gamma_1$ is estimated to be 0.155, significant at level 5\%. However, the model fails to capture the buffering effect of self-efficacy in this relationship and $\gamma_1$ is a constant. The effect of depressed affect on withdrawal intentions $\alpha_0$ is estimated to be 1.000, significant at the 0.05 level. Thus, the indirect effect of economic stress on entrepreneurs' intentions to withdraw $\beta=\gamma_1\alpha_0=0.155$, with a 95\% confidence interval $(0.054, 0.256)$. It is noteworthy that when the dynamic pattern of covariates is ignored, the effect of gender on depressed affect is no longer significant. One possible reason is that positive and negative effects of gender at different self-efficacy levels cancel out with each other, leading to an insignificant estimate.

\section{Conclusion}\label{sec5}

This study proposes a generalized varying coefficient mediation model (GVCMM) that can be used to examine the indirect effect of economic stress on entrepreneurial withdrawals via depressed affect. Different from standard mediation models, GVCMM allows indirect effect and direct effect to be represented as smooth functions of a continuous variable, and accommodates various types of responses. Estimation procedures and asymptotic properties are established. Simulation studies are conducted to verify the finite sample performance of the proposed methodology. The real data analysis for behavioral economics shows that self-efficacy serves to buffer the deleterious impact of economic stress, which in turn reduces business owners' withdrawal intentions. In addition to behavioral economics, GVCMM can also be applied to studies in various fields.

\section{Technical Proofs of the Main Results}
\subsection{Proofs of Theorem \ref{GVCMMalpha1_u} and Lemma \ref{GVCMMlem_alpha0_u}}\label{secA1}
By asymptotic properties of generalized linear model, we derive that
\begin{equation*}
\sqrt{n}\left\{\left(
   \begin{array}{c}
      \hat{\bm{\alpha}}_0^*\\
      \hat{\bm{\alpha}}_1^*
   \end{array}
 \right)-
 \left(
   \begin{array}{c}
      \bm{\alpha}_0^*\\
      \bm{\alpha}_1^*
   \end{array}
 \right)\right\}\overset{D}{\rightarrow}
 N\left\{\bm{0},\phi \bm{I}_1^{-1}(\hat{\bm{\alpha}}^{*T}, \hat{\bm{\alpha}}_1^{*T})\right\}\overset{D}{=}N(\bm{0},\phi\bm{\Sigma^{*-1}}),
\end{equation*}
where
$\bm{\Sigma}^*=\textbf{I}_1(\hat{\bm{\alpha}}^{*T}, \hat{\bm{\alpha}}_1^{*T})=\mathrm{E}\{\bm{z}^{*}_i\bm{z}^{*T}_i b^{''}(\bm{\alpha}_0^{*T}\bm{m}_i^{*} +\bm{\alpha}_1^{*T}\bm{x}_i^{*})\},$  $\bm{z}^*_i=(\bm{m}_i^{*T},\bm{x}_i^{*T})^T.$ To be specific,
\begin{align*}
\bm{\Sigma}^*=
    \left(\begin{array}{cc}
        \bm{\Sigma}_{m^*m^*} & \bm{\Sigma}_{m^*x^*} \\
        \bm{\Sigma}_{x^*m^*} & \bm{\Sigma}_{x^*x^*}
    \end{array}\right)
\end{align*}
and
\begin{align*}
    \bm{\Sigma}_{m^*m^*}&=\mathrm{E}\{\bm{m}_i^*\bm{m}_i^{*T}b^{''}(\bm{\alpha}_0^{*T}\bm{m}_i^*+\bm{\alpha}_1^{*T}\bm{x}_i^*)\},\\
    \bm{\Sigma}_{m^*x^*}&=\mathrm{E}\{\bm{m}_i^*\bm{x}_i^{*T}b^{''}(\bm{\alpha}_0^{*T}\bm{m}_i^*+\bm{\alpha}_1^{*T}\bm{x}_i^*)\},\\
    \bm{\Sigma}_{x^*m^*}&=\mathrm{E}\{\bm{x}_i^*\bm{m}_i^{*T}b^{''}(\bm{\alpha}_0^{*T}\bm{m}_i^*+\bm{\alpha}_1^{*T}\bm{x}_i^*)\},\\
    \bm{\Sigma}_{x^*x^*}&=\mathrm{E}\{\bm{x}_i^*\bm{x}_i^{*T}b^{''}(\bm{\alpha}_0^{*T}\bm{m}_i^*+\bm{\alpha}_1^{*T}\bm{x}_i^*)\}.
\end{align*}

Therefore, the asymptotic distributions of $\hat{\bm{\alpha}}_0^*$ and $\hat{\bm{\alpha}}_1^*$ can be derived as
$$\sqrt{n}(\hat{\bm{\alpha}}_0^*-\bm{\alpha}_0^*)\overset{D}{\rightarrow}N(\bm{0},\phi\bm{\Sigma}_{m^*m^*\mid x^*}^{-1}),$$
$$\sqrt{n}(\hat{\bm{\alpha}}_1^*-\bm{\alpha}_1^*)\overset{D}{\rightarrow}N(\bm{0},\phi\bm{\Sigma}_{x^*x^*\mid m^*}^{-1}),$$
where
\begin{align}
    \label{Sigm*m*.x*}\bm{\Sigma}_{m^*m^*\mid x^*}&=\bm{\Sigma}_{m^*m^*}-\bm{\Sigma}_{m^*x^*}\bm{\Sigma}_{x^*x^*}^{-1}\bm{\Sigma}_{x^*m^*},\\
    \label{Sigx*x*.m*}\bm{\Sigma}_{x^*x^*\mid m^*}&=\bm{\Sigma}_{x^*x^*}-\bm{\Sigma}_{x^*m^*}\bm{\Sigma}_{m^*m^*}^{-1}\bm{\Sigma}_{m^*x^*}.
\end{align}

Thus, for any $q\times qK_1$ non-zero constant matrix $G_1$ and $p\times pK_0$ non-zero constant matrix $G_0$,
$$\{\Var(\bm{G}_1\hat{\bm{\alpha}}_1^*)\}^{-1/2}\bm{G}_1(\hat{\bm{\alpha}}_1^*-\bm{\alpha}_1^*)\overD N_q(\0,\bm{I}),$$
$$\{\Var(\bm{G}_0\hat{\bm{\alpha}}_0^*)\}^{-1/2}\bm{G}_0(\hat{\bm{\alpha}}_0^*-\bm{\alpha}_0^*)\overD N_p(\0,\bm{I}),$$
where $\Var(\bm{G}_1\hat{\bm{\alpha}}_1^*) = \frac{1}{n}\phi\bm{G}_1\bm{\Sigma}_{x^*x^*\mid m^*}^{-1}\bm{G}_1^T$ and $\Var(\bm{G}_0\hat{\bm{\alpha}}_0^*) = \frac{1}{n}\phi\bm{G}_0\bm{\Sigma}_{m^*m^*\mid x^*}^{-1}\bm{G}_0^T$.

\subsection{Proof of Lemma \ref{GVCMMlem_Gamma_u}}\label{secA2}
To get the estimator of $\bm{C}$, we need to minimize $\ell_2=\left\Vert\bm{M}-\bm{X}^m\bm{C}\right\Vert^2$ with respect to $\bm{C}$. Suppose that $\bm{X}^m\bm{X}^{mT}$ is invertible, then $\ell_2(\bm{C})$ has a unique minimizer and
\begin{align*}
\hat{\bm{C}}=&(\bm{X}^{mT}\bm{X}^m)^{-1}\bm{X}^{mT}\bm{M}\\
=&(\bm{X}^{mT}\bm{X}^m)^{-1}\bm{X}^{mT}(\bm{X}^m\bm{C}+\bm{E})\\
=&\bm{C}+(\bm{X}^{mT}\bm{X}^m)^{-1}\bm{X}^{mT}\bm{E}.
\end{align*}

To derive the asymptotic distribution of $\hat{\bm{C}}$, let us define the variance of a matrix at first. Let $\bm{B}$ be a matrix of dimension $N\times p$, where $\bm{B}=(\bm{b}_1,\cdots, \bm{b}_p)$. Then $\mathrm{vec}(\bm{B})=(\bm{b}_1^T,\cdots, \bm{b}_p^T)^T$, and $\Var(\bm{B})=\Var\{\mathrm{vec}(\bm{B}^T)\}$.

By Weak Law of Large Numbers and Central Limit Theorem, for any $q\times qK_m$ constant matrix $\bm{G}_2$, $\bm{G}_2 (\hat{\bm{C}}-\bm{C})=\frac{1}{n}\bm{G}_2(\frac{1}{n}\bm{X}^{mT}\bm{X}^m)^{-1}\bm{X}^{mT}\bm{E}$ follows normality with mean zero. Since
\begin{align}
&\mathrm{Var}\{\bm{G}_2 (\hat{\bm{C}}-\bm{C})\}-\frac{1}{n}(\bm{G}_2\bm{\Sigma}_{x^mx^m}^{-1}\bm{C}_n^{'T})\otimes \bm{\Sigma}_\epsilon\\
\overset{p}{\rightarrow}&\frac{1}{n^2}\mathrm{Var}\{\mathrm{vec}(\bm{I}_q\bm{E}^T\bm{X}^m\bm{\Sigma}_{x^mx^m}^{-1}\bm{G}_2^{T})\} - \frac{1}{n}(\bm{G}_2\bm{\Sigma}_{x^mx^m}^{-1}\bm{C}_n^{'T})\otimes \bm{\Sigma}_\epsilon\nonumber\\
=&\frac{1}{n^2}\mathrm{Var}\{(\bm{G}_2\bm{\Sigma}_{x^mx^m}^{-1}\bm{X}^{mT}\otimes \bm{I}_q)\mathrm{vec}(\bm{E}^T)\} - \frac{1}{n}(\bm{G}_2\bm{\Sigma}_{x^mx^m}^{-1}\bm{C}_n^{'T})\otimes \bm{\Sigma}_\epsilon\nonumber\\
\overset{p}{\rightarrow}&0,\nonumber
\end{align}
we obtain that $$\{\Var(\bm{G}_2\hat{\bm{C}})\}^{-1/2}\bm{G}_2(\hat{\bm{C}}-\bm{C})\overD N(\bm{0},\bm{I}),$$ where $\Var(\bm{G}_2\hat{\bm{C}}) = \frac{1}{n}(\bm{G}_2\bm{\Sigma}_{x^mx^m}^{-1}\bm{G}_2^{T})^{-1}\otimes \bm{\Sigma}_\epsilon$.

\subsection{Proof of Theorem \ref{GVCMMrates}}\label{secA3}
Under Conditions $\bm{C1-C2}$, the following properties of B-splines are valid.

\begin{itemize}
    \item(\citealt{de1978practical} and \citealt{huang2004polynomial}) For $k=1,\cdots, K_0$, $b_{0k}(\cdot)$ are the B-spline basis functions that span $\mathbb{G}_{\alpha_{0j}}$, for any $j=1,\cdots, p$, then $b_{0k}(u)\geq 0$ and $\sum_{k=1}^{K_0}b_{0k}(u)=1$, $u\in\mathcal{U}.$ In addition, there exist positive constants $C_2$ and $C_3$ such that
    $$\frac{C_2}{K_0}\sum_{k}\gamma_{0k}^2\leq \int_\mathcal{U}\left\{\sum_k \gamma_{0k}b_{0k}(u)\right\}^2 f_U(u) du\leq \frac{C_3}{K_0}\sum_k\gamma_{0k}^2,$$ $$\gamma_{0k}\in\mathbb{R}, k=1,\cdots, K_0.$$
    The mentioned property also holds for $b_{1k}(\cdot)$, $k=1,\cdots, K_1,$ and $B(\cdot)$, $k=1,\cdots, K_m$, correspondingly.
\end{itemize}

Based on Conditions \textbf{C1-C8}, as explained in \cite{huang2004polynomial}, we can obtain that for any $j=1,\cdots, p$ and $l=1,\cdots, q$, $\left\Vert \hat{\Gamma}_{lj}-\Gamma_{lj} \right\Vert^2_{L_2}=O_p(1/n+K_m/n)=O_p(1/n)$.

Then, we want to get the convergence rate of $\hat{\bm{\alpha}}_0$. Since
\begin{align*}
    \left\Vert \hat{\alpha}_{0j}-\alpha_{0j}\right\Vert_{L_2}^2&= \int_\mathcal{U}\mid \hat{\alpha}_{0j}(\cdot)-\alpha_{0j}(\cdot)\mid ^2f_U(u)du\\
    &=\int_\mathcal{U}\mid (\hat{\bm{a}}_{0j}-\bm{a}_{0j})^T\bm{b}_0(u)\mid ^2f_U(u)du\\
    &=(\hat{\bm{a}}_{0j}-\bm{a}_{0j})^T\int_\mathcal{U}\bm{b}_0(u)\bm{b}_0(u)^Tf_U(u)du(\hat{\bm{a}}_{0j}-\bm{a}_{0j})\\
    &\leq \lambda_{max}(\bm{A})\left\Vert (\hat{\bm{a}}_{0j}-\bm{a}_{0j})\right\Vert^2
\end{align*}
where the matrix $\bm{A}=\int_\mathcal{U}\bm{b}_0(u)\bm{b}_0(u)^Tf_U(u)du$, and $\lambda_{max}(\cdot)$ denotes the largest eigenvalue of a matrix. Since the matrix $\bm{A}$ is a positive square matrix, by Collatz-Wielandt formula, the largest eigenvalue of it is bounded above by the maximal row sum of $\bm{A}$, which is bounded away from infinity by the properties of B splines.
As a result, once we know $\left\Vert\hat{\bm{a}}_{0j}-\bm{a}_{0j}\right\Vert^2=O_p(1/n)$, we can derive that $\left\Vert \hat{\alpha}_{0j}-\alpha_{0j}\right\Vert_{L_2}^2=O_p(1/n).$

We want to show that for any given $\epsilon>0$, there exists a large constant C such that
\begin{align}
    \label{prob_inequality}P\{\sup_{\left\Vert\bm{v}\right\Vert=C}L(\bm{\theta}^{*}+r_n\bm{v})<L(\bm{\theta}^{*})\}\geq 1-\epsilon
\end{align}
where $\bm{\theta}^{*}=(\bm{\alpha}_0^{*T},\bm{\alpha}_1^{*T})^T$ and $L(\bm{\theta}^{*})=L(\bm{\alpha}_0^{*}, \bm{\alpha}_1^{*})$ is the log-likelihood function, $\bm{v}\in\mathbb{R}^{q+p}$, and $r_n=1/\sqrt{n}.$ This implies with probability at least $1-\epsilon$, there exists a maximum in the ball $\{\bm{\theta}^{*}+r_n\bm{v}:\left\Vert\bm{v}\right\Vert<C\}.$ Hence, there exists a local maximizer $\hat{\bm{\theta}}^*$ such that $\left\Vert \hat{\bm{\alpha}}_0^{*}-\bm{\alpha}_0^{*}\right\Vert=O_p(r_n)$ and  $\left\Vert \hat{\bm{a}}_{0j}^{*}-\bm{a}_{0j}^{*}\right\Vert=O_p(r_n)$, for all $j=1,\cdots, p.$

By Taylor expansion,
\begin{align*}
    &L(\bm{\theta}^*+r_n\bm{v})-L(\bm{\theta}^*)\\
    \leq & r_nL^{'}(\bm{\alpha}_0^{*},\bm{\alpha}_1^{*})^T\bm{v}-\frac{1}{2}\bm{v}^TI_1(\bm{\alpha}_0^{*},\bm{\alpha}_1^{*})\bm{v}nr_n^2\{1+o_p(1)\}\\
    \leq & r_nL^{'}(\bm{\alpha}_0^{*},\bm{\alpha}_1^{*})^T\bm{v}-\frac{1}{2}\lambda_{min}\{I_1(\bm{\alpha}_0^{*},\bm{\alpha}_1^{*})\}\left\Vert\bm{v}\right\Vert^2nr_n^2\{1+o_p(1)\}.
\end{align*}
Note that $n^{-1/2}L^{'}(\bm{\alpha}_0^{*},\bm{\alpha}_1^{*})=O_p(1),$ $L^{'}(\bm{\alpha}_0^{*},\bm{\alpha}_1^{*})=O_p(n^{1/2})$, and according to Condition \textbf{C8}, by choosing a sufficiently large $C$, the second term dominates the first term uniformly in $\left\Vert\bm{v}\right\Vert=C.$ So the inequality (\ref{prob_inequality}) has been proved. As a result, $\left\Vert \hat{\alpha}_{0j}-\alpha_{0j}\right\Vert_{L_2}^2=O_p(1/n).$

Since
\begin{align*}
    \left\Vert \hat{\beta}_l-\beta_l \right\Vert^2_{L_2}&=\int_{u\in\mathcal{U}}\mid \hat{\beta}_l(u)-\beta_l(u)\mid ^2 f_U(u)du\\
    &=\int_{u\in\mathcal{U}}\mid \sum_{j=1}^p\{\hat{\Gamma}_{lj}(u)\hat{\alpha}_{0j}(u)-\Gamma_{lj}(u)\alpha_{0j}(u)\}\mid ^2f_U(u) du\\
    &\leq p\cdot \sum_{j=1}^p \int_{u\in\mathcal{U}}\mid \hat{\Gamma}_{lj}(u)\hat{\alpha}_{0j}(u)-\Gamma_{lj}(u)\alpha_{0j}(u)\mid ^2 f_U(u)du
\end{align*}
by Cauchy Inequality, and
\begin{align*}
    &\int_{u\in\mathcal{U}}\mid \hat{\Gamma}_{lj}(u)\hat{\alpha}_{0j}(u)-\Gamma_{lj}(u)\alpha_{0j}(u)\mid^2 f_U(u)du\\
    =&\int_{u\in\mathcal{U}}\mid \hat{\Gamma}_{lj}(u)\{\hat{\alpha}_{0j}(u)-\alpha_{0j}(u)\} + \alpha_{0j}(u)\{\hat{\Gamma}_{lj}(u)-\Gamma_{lj}(u)\}\mid ^2 f_U(u)du\\
    \leq & 2\int_{u\in\mathcal{U}}\mid \hat{\Gamma}_{lj}(u)\{\hat{\alpha}_{0j}(u)-\alpha_{0j}(u)\}\mid ^2 f_U(u)du \\
              &+ 2\int_{u\in\mathcal{U}}\mid \alpha_{0j}(u)\{\hat{\Gamma}_{lj}(u)-\Gamma_{lj}(u)\}\mid ^2 f_U(u)du\\
    \leq & 2M_{\hat{\Gamma}}^2\int_{u\in\mathcal{U}}\mid \hat{\alpha}_{0j}(u)-\alpha_{0j}(u)\mid ^2 f_U(u)du + 2M_0^2\int_{u\in\mathcal{U}}\mid \hat{\Gamma}_{lj}(u)-\Gamma_{lj}(u)\mid ^2 f_U(u)du\\
    =&2M_{\hat{\Gamma}}^2\left\Vert \hat{\alpha}_{0j}-\alpha_{0j} \right\Vert^2_{L_2} + 2M_0^2\left\Vert \hat{\Gamma}_{lj}-\Gamma_{lj} \right\Vert^2_{L_2}
\end{align*}
by Condition \textbf{C5}, we can derive that
\begin{align*}
\int_{u\in\mathcal{U}}\mid \hat{\Gamma}_{lj}(u)\hat{\alpha}_{0j}(u)-\Gamma_{lj}(u)\alpha_{0j}(u)\mid ^2 f_U(u)du=O_p(1/n).
\end{align*}
It's obvious that $\left\Vert\hat{\beta}_l-\beta_l\right\Vert^2_{L_2}=O_p(1/n).$

\subsection{Proof of Theorem \ref{GVCMMasy.}}\label{secA4}
The asymptotic distribution of $\hat{\bm{\beta}}(u)$ can be obtained as follows. For any $\bm{a}_{q\times 1}$, we consider the Cramer's Device. According to the result in Lemma \ref{GVCMMlem_Gamma_u},
$$\Var\{\mathrm{vec}(\hat{\bm{\Gamma}}(u)^T-\bm{\Gamma}(u)^T)\}=\frac{1}{n}\big[\{\bm{I}_q\otimes \bm{B}(u)^T\}\bm{\Sigma}_{x^mx^m}^{-1}\{\bm{I}_q\otimes \bm{B}(u)\}\big]\otimes\bm{\Sigma}_\epsilon.$$
Since
$\mathrm{vec}(\bm{I}_p\hat{\bm{\Gamma}}(u)^T\bm{a})=(\bm{a}^T\otimes \bm{I}_p)\mathrm{vec}(\hat{\bm{\Gamma}}(u)^T)$
and
$\mathrm{vec}(\hat{\bm{\Gamma}}(u)^T\bm{a}-\bm{\Gamma}(u)^T\bm{a})=(\bm{a}^T\otimes \bm{I}_p)\{\mathrm{vec}(\hat{\bm{\Gamma}}(u)^T)-\mathrm{vec}(\bm{\Gamma}(u)^T)\},$
we can get that
\begin{align*}
&\Var\{\mathrm{vec}(\hat{\bm{\Gamma}}(u)^T\bm{a}-\bm{\Gamma}(u)^T\bm{a})\}\\=&\Var[(\bm{a}^T\otimes \bm{I}_p)\{\mathrm{vec}(\hat{\bm{\Gamma}}(u)^T-\bm{\Gamma}(u)^T)\}]\nonumber\\
=&(\bm{a}^T\otimes\bm{I}_p)\bigg(\frac{1}{n}\big[\{\bm{I}_q\otimes \bm{B}(u)^T\}\bm{\Sigma}_{x^mx^m}^{-1}\{\bm{I}_q\otimes \bm{B}(u)\}\big]\otimes\bm{\Sigma}_\epsilon\bigg)(\bm{a}\otimes\bm{I}_p)\nonumber\\
=&\frac{1}{n}\big[\bm{a}^T\{\bm{I}_q\otimes \bm{B}(u)^T\}\bm{\Sigma}_{x^mx^m}^{-1}\{\bm{I}_q\otimes \bm{B}(u)\}\bm{a}\big]\otimes \bm{\Sigma}_\epsilon.\nonumber\\
\end{align*}
If we let $\bm{\theta}(u)=\bm{\Gamma}(u)^T \bm{a}$ and $\hat{\bm{\theta}}(u)=\hat{\bm{\Gamma}}(u)^T\bm{a}$, then
$$[\Var\{\hat{\bm{\theta}}(u)\}]^{-1/2}(\hat{\bm{\theta}}(u)^T-\bm{\theta}(u)^T)\overset{D}{\rightarrow}N(\bm{0},\bm{I}),$$
where $\Var\{\hat{\bm{\theta}}(u)\} = \frac{1}{n}\big[\bm{a}^T\{\bm{I}_q\otimes \bm{B}(u)^T\}\bm{\Sigma}_{x^mx^m}^{-1}\{\bm{I}_q\otimes \bm{B}(u)\}\bm{a}\big]\bm{\Sigma}_\epsilon$.

Now we want to get the asymptotic distribution of $\hbe(u)=\hGam(u) \halp(u)$ by delta method.  Let $\bm{\Sigma}_{\theta(u)} = \Var\{\hat{\bm{\theta}}(u)\}$ and $\bm{\Sigma}_{\bm{\alpha}_0(u)} = \Var\{\hat{\bm{\alpha}}_0(u)\}$ for simplicity.

Since $\halp(u)$ and $\hGam(u)$ are independent and
\begin{equation*}
\left(
    \begin{array}{cc}
         \bm{\Sigma}_{\theta(u)} &\bm{0}\\
         \bm{0}& \bm{\Sigma}_{\alpha_0(u)}\\
    \end{array}
  \right)^{-1/2}
\left\{\left(
   \begin{array}{c}
      \hat{\bm{\theta}}(u)^T\\
      \hat{\bm{\alpha}}_0(u)
   \end{array}
 \right)-
 \left(
   \begin{array}{c}
      \bm{\theta}(u)^T\\
      \bm{\alpha}_0(u)
   \end{array}
 \right)\right\}\overset{D}{\rightarrow}
 N\left(\bm{0},\bm{I}\right),
\end{equation*}
we have $$\{\bm{\alpha}_0(u)^T\bm{\Sigma}_{\theta(u)}\bm{\alpha}_0(u)+\bm{\theta}(u)^T\bm{\Sigma}_{\alpha_0(u)}\bm{\theta}(u)\}^{-1/2}\{\hat{\bm{\theta}}(u)^T\hat{\bm{\alpha}}_0(u)-\bm{\theta}(u)^T\bm{\alpha}_0(u)\}\overset{D}{\rightarrow}N(\bm{0},\bm{I}),$$
where $$\bm{\alpha}_0(u)^T\bm{\Sigma}_{\theta(u)} \bm{\alpha}_0(u) = \frac{1}{n}\bm{\alpha}_0(u)^T\Sigep \bm{\alpha}_0(u) \bm{a}^T\{\bm{I}_q\otimes \bm{B}(u)^T\}\bm{\Sigma}_{x^mx^m}^{-1}\{\bm{I}_q\otimes \bm{B}(u)\}\bm{a}$$ and $$\bm{\theta}(u)^T\bm{\Sigma}_{\alpha_0(u)}\bm{\theta}(u)=\frac{1}{n}\phi\bm{a}^T\bm{\Gamma}(u)\{\bm{I}_p\otimes \bm{b}_0(u)^T\}\bm{\Sigma}_{m*m*\mid x*}^{-1}\{\bm{I}_p\otimes \bm{b}_0(u)\}\bm{\Gamma}(u)^T\bm{a}.$$
Thus,
\begin{align*}
&\bm{\alpha}_0(u)^T\bm{\Sigma}_{\theta(u)}\bm{\alpha}_0(u)+\bm{\theta}(u)^T\bm{\Sigma}_{\alpha_0(u)}\bm{\theta}(u)\\
=&\frac{1}{n}\bm{a}^T\big[\bm{\alpha}_0(u)^T\Sigep \bm{\alpha}_0(u) \{\bm{I}_q\otimes \bm{B}(u)^T\}\bm{\Sigma}_{x^mx^m}^{-1}\{\bm{I}_q\otimes \bm{B}(u)\}\\
    &+\phi\bm{\Gamma}(u)\{\bm{I}_p\otimes \bm{b}_0(u)^T\}\bm{\Sigma}_{m*m*\mid x*}^{-1}\{\bm{I}_p\otimes \bm{b}_0(u)\}\bm{\Gamma}(u)^T\big]\bm{a}.
\end{align*}
Since $\{\hat{\bm{\theta}}(u)^T\halp(u)-\bm{\theta}(u)^T\bm{\alpha}_0(u)\}=\bm{a}^T\{\hGam(u)\halp(u)-\bm{\Gamma}(u)\bm{\alpha}_0(u)\}$, we derive that
\begin{align*}
&[\Var\{\hat{\bm{\beta}}(u)\}]^{-1/2}\{\hat{\bm{\Gamma}}(u)\hat{\bm{\alpha}}_0(u)-\bm{\Gamma}(u)\bm{\alpha}_0(u)\}\\
=&[\Var\{\hat{\bm{\beta}}(u)\}]^{-1/2}\{\hat{\bm{\beta}}(u)-\bm{\beta}(u)\}\\
\overset{D}{\rightarrow}&N\{\bm{0},\bm{I}\},
\end{align*}
where $\Var\{\hat{\bm{\beta}}(u)\}=\frac{1}{n}[\bm{\alpha}_0(u)^T\Sigep \bm{\alpha}_0(u) \{\bm{I}_q\otimes \bm{B}(u)^T\}\bm{\Sigma}_{x^mx^m}^{-1}\{\bm{I}_q\otimes \bm{B}(u)\}+\phi\bm{\Gamma}(u)\{\bm{I}_p\otimes \bm{b}_0(u)^T\}\bm{\Sigma}_{m*m*\mid x*}^{-1}\{\bm{I}_p\otimes \bm{b}_0(u)\}\bm{\Gamma}(u)^T].$


\backmatter

\bmhead{Acknowledgments}



Liu's research was supported by National Natural Science Foundation of China 11771361 and 11871409.





\bibliography{GVCMM}

\end{document}